\documentclass[a4paper]{article}
\usepackage{lineno,hyperref}
\usepackage{amsmath,amssymb}
\usepackage{cases}
\usepackage{color}
\usepackage{graphicx}
\usepackage{comment}
\begin{document}

\title{Boundary integral equation methods for the two dimensional wave equation in time domain revisited}
\author{
Mio Fukuhara$^{1)}$\thanks{fukuhara@acs.i.kyoto-u.ac.jp} \and
Ryota Misawa$^{2)}$ \thanks{misawa.ryota.8c@kyoto-u.ac.jp} \and
Kazuki Niino$^{1)}$ \thanks{niino@acs.i.kyoto-u.ac.jp} \and
Naoshi Nishimura$^{1)}$\thanks{corresponding author, nchml@i.kyoto-u.ac.jp}}
\date{}
\maketitle

$1)$ Graduate School of Informatics, Kyoto University, Kyoto 606-8501, Japan

$2)$ Graduate School of Engineering, Kyoto University, Kyoto 615-8540, Japan

\section*{abstract}
  This study considers the stability of time domain BIEMs for the wave
  equation in 2D. We show that the stability of time domain BIEMs is
  reduced to a nonlinear eigenvalue problem related to frequency
  domain integral equations. We propose to solve this non-linear
  eigenvalue problem numerically with the Sakurai-Sugiura
  method. After validating this approach numerically in the exterior
  Dirichlet problem, we proceed to transmission problems in which we
  find that some time domain counterparts of ``resonance-free''
  integral equations in frequency domain lead to instability. We
  finally show how to reformulate these equations to obtain stable
  numerical schemes. 

\section*{keywords}
stability, time domain, BIEM, transmission problems, eigenvalue problems

\section{Introduction}
Boundary integral equation methods (BIEMs) are often said to be
advantageous in wave problems because they can be applied to
scattering problems easily.  It is certainly true that BIEMs in
frequency domain are easy to use, but the same does not necessarily
apply to time domain methods.  As a matter of fact, BIEMs for the wave
equation in time domain have a long standing stability problem and
there have been many efforts to stabilise BIEMs for wave equations.
For example, Ha Duong and his colleagues (e.g., \cite{bamberger})
showed the stability of some time domain BIEMs in 3D based on
space-time variational formulations. Their argument depends on the
energy conservation which is why their variational formulation
includes time derivatives (e.g., the time derivative of single layer
potential). Aimi et al.\ presented some numerical results in 2D using
time or space differentiated integral equations and a fully
variational approach.  Abboud et al.\cite{abboud} considered a
coupling of fully variational BIEMs with discontinuous Galerkin
methods.  Unfortunately, however, implementing computational codes for
the full space-time variational formulation is not very easy.  Coding
becomes easier if one uses variational approaches only spatially and
use collocation in time. Van 't Wout et al.\cite{vantwout} have shown
a way to find a stable time-collocated variational approach based on
fully variational methods. In spite of these efforts, the standard
collocation approaches remain the preferred choice in engineering,
although known mathematical stability results in collocation are
rather limited (see Davies and Dancan\cite{davies} for example).
Various numerical stabilisation techniques for collocation have been
proposed, from which we cite just a few relatively new
ones\cite{parot1,parot2,jang1,jang2,pak,bai}, referring the reader to
the lists of references of these papers for further literatures.  Some
other investigations take viewpoints similar to ours in that they seek
stabilisation based on the choices of integral equations.  For
example, the use of time differentiated integral equations has been
advocated by several authors\cite{mansur, panag}. Ergin et
al.\cite{ergin} proposed to use the Burton-Miller (BM) integral
equation to achieve stability guided by an observation that the
instability of BIEMs for scattering problems is related to fictitious
eigenfrequencies (internal resonance). Chappell et
al.\cite{chappell,chappell2} gave further insight as well as the
implementation details of the BM formulation. This formulation has
been utilised recently in practical applications\cite{zhang}.  Finally we
mention recent developments of CQM by
Lubich\cite{lubich,bb:sayas,banjai} which is a stable method of
computing convolutions. CQM has been applied successfully to
engineering applications (e.g., Schanz et al. \cite{schanz}). However,
implementing CQM is still not as simple as the standard collocation
methods, which is the reason we consider the conventional approach in
this paper.

The above brief review of the works on the stability of time domain
BIEMs for the wave equation covers just a small part of what have been
done so far. Indeed, the cause of the instability is now fairly well
understood in connection with the spectra of the integral operators
and the error introduced by discretisation (e.g.,
\cite{ergin,chappell,wang}), particularly in exterior problems. In
spite of these efforts by predecessors, however, there seem to exist
no definite and simple criteria of stability for the collocation
methods.  One still needs to carry out a quantitative assessment
numerically in order to see if a particular scheme is stable or not.
A standard method to check the stability of collocation BIEMs in time
domain is to compute characteristic roots by solving a polynomial
eigenvalue problem (see \eqref{algeh}) after reducing it to an
equivalent linear eigenvalue problem for the companion matrix (See,
e.g., Walker et al.\cite{walker}). This method is particularly
effective in 3D where the fundamental solution has a finite ``tail''
(i.e., it vanishes after a finite time). However, this approach needs
linear eigensolvers for sparse, but large, matrices. One may possibly
solve polynomial eigenvalue problems directly to reduce the size of
the matrix, but this will lead to a non-linear eigenvalue
problem. Fortunately, recent developments of eigensolvers based on
contour integrals such as the Sakurai-Sugiura method (SSM)
\cite{bb:ssm} made the solution of non-linear eigenvalue problems
feasible. In 2D problems, however, the same approach is not very
practical because the fundamental solution in 2D is very slow to decay
in time.  In this paper we propose to resolve this difficulty by
carrying out the required stability analysis in frequency domain.
Namely, we convert the stability analysis for BIEMs in two dimensional
wave equation to a non-linear eigenvalue problem similar to those for
the Helmholtz equation and solve it with SSM using techniques proposed
in Misawa et al.\cite{misawa1,misawa2}.  This approach has an
additional benefit of making the relation between eigenvalues of the
approximated integral operators in frequency domain and the stability
clearer, thus providing new intuition into the subject. Using the
proposed technique, we investigate stability of various time domain
integral equations for transmission problems, which have not been
investigated very much so far.

As a basic study in this subject, however, this paper considers very
simple problems only. Namely, we restrict our attention mainly to
exterior Dirichlet problems and transmission problems for domains
bounded by a circle. We first present a stability analysis for the
exterior Dirichlet problem which uses frequency domain tools.  The
question of stability is then reduced to the computation of the
characteristic roots which are eigenvalues of a certain non-linear
equation. After solving this eigenvalue problem with SSM ignoring the
effect of the spatial discretisation, we identify potentials which
yield stable numerical schemes with piecewise linear time basis
functions. We then proceed to transmission problems in which we show
that even the time domain counterparts of ``resonance free'' BIEMs may
lead to instability. We then modify these integral equations by using
only the ``stable potentials'' and show by numerical experiments that
these formulations do lead to stability in time domain. This
conclusion is supported by the stability analysis using SSM. After
examining the influence of the spatial discretisation on the
characteristic roots, we present numerical examples of transmission
problems for non-circular domains solved with the modified integral
equations, which turn out to be stable.

\section{Exterior Dirichlet problems}
\subsection{Formulation}
\label{sec:edp}
Let $D_2 \subset \mathbb{R}^2$ be a bounded domain whose boundary $\Gamma=\partial D_2$ is smooth and let $D_1$ be the exterior of $D_2$, i.e., $D_1=\mathbb{R}^2\setminus \bar D_2$.
Also, let $n$ be the unit normal vector on $\Gamma$ directed towards $D_1$. 
We are interested in the following initial- boundary value problem
(Dirichlet problem):

Find $u$ which satisfies the two dimensional wave equation in $D_1$:
\begin{align}
\label{wave}
\Delta u	-\frac{1}{c_1^2}\frac{\partial^2 u}{{\partial t}^2} = 0 \ \ \ \mathrm{in} \  D_1 \times (t>0)
\end{align}
the homogeneous Dirichlet boundary condition on $\Gamma$ for $u$:
\begin{align*}
u= 0 \ \ \ \mathrm{on} \  \Gamma \times (t>0)
\end{align*}
the homogeneous initial conditions in $D_1$:
\begin{align}
\label{eq:init}
u^{\mathrm{sca}}|_{t=0}=\frac{\partial u^{\mathrm{sca}}}{\partial t}|_{t=0}= 0, \ \ \ \mathrm{in} \  D_1
\end{align}
and the radiation condition for the scattered wave $u^{\mathrm{sca}}=u-u^{\mathrm{inc}}$ in $D_1$, where $c_1$ is the wave speed in $D_1$ which is written as $c_1=\sqrt{\frac{s_1}{\rho_1}}$, $s_1$ and $\rho_1$ are the shear modulus and density in $D_1$ and $u^{\mathrm{inc}}$ is the incident wave which satisfies \eqref{wave} in the whole space-time, respectively.

\subsection{Boundary integral equations}
The solution to the above initial- boundary value problem can be written as
\begin{align*}
u(x,t)&=u^{\mathrm{inc}}(x,t)-S^1q(x,t)\qquad (x,t) \in D_1\times (t>0)
\end{align*}
if the function $q(x,t)$ on $\Gamma \times (t>0)$ is chosen such that
\begin{align}
  \label{solin}
0&=u^{\mathrm{inc}}(x,t)-S^1q(x,t)\qquad (x,t) \in D_2\times (t>0)
\end{align}
is satisfied, where
$S^{\nu}$ ($\nu=1$ or 2. $\nu=1$ in the
present context) stands for the single layer potential defined by
\begin{align*}
S^\nu q(x,t) := \int^t_0 \int_\Gamma G^\nu (x-y,t-s)q(y,s) \ dS_yds
\end{align*}
and  $G^\nu (x,t)$ is the fundamental solution of the wave equation given by:
\begin{align}
\label{eq:fs}
  G^\nu (x,t)=\frac {c_\nu}{2 \pi \sqrt{(c_\nu t)^2 - |x|^2}_+}=
\begin{cases}
\frac {c_\nu}{2 \pi\sqrt{(c_\nu t)^2 - |x|^2}}& c_\nu t>|x|\\
0&\mbox{otherwise.}
\end{cases}
\end{align}
For later convenience, we also introduce the normal derivative of the
single layer potential $D^{T\nu}$, the double layer potential $D^\nu$
and its normal derivative $N^{\nu}$ defined by:
\begin{align*}
D^{T\nu} q(x,t) :=\int^t_0 \int_\Gamma \frac{\partial G^\nu}{\partial n_x} (x-y,t-s)q(y,s) \ dS_yds,\\
D^\nu q(x,t) :=\int^t_0 \int_\Gamma \frac{\partial G^\nu}{\partial n_y} (x-y,t-s)q(y,s) \ dS_yds, \\
N^{\nu} q(x,t) :=\int^t_0 \int_\Gamma \frac{\partial^2 G^\nu}{\partial n_x\partial n_y} (x-y,t-s)q(y,s) \ dS_yds.
\end{align*}
By these notations for potentials we indicate functions defined in
$x\in \mathbb{R}^2\setminus \Gamma$ in this paper. Their boundary
traces from the exterior (interior) are indicated by superposed $+$
($-$). When the exterior and interior traces
coincide, however, we denote them by the same letter without superposed $\pm$.
This applies to $S^{\nu}$, $(d/dt) S^{\nu}=\dot S^{\nu}$ and $N^{\nu}$, but we need to evaluate
the boundary integral in  $N^{\nu}$ in the sense of the finite part then.

The function $q(x,t)$ coincides with the (exterior trace of)
normal derivative of $u$ on $\Gamma$ if \eqref{solin} is satisfied.
The condition in \eqref{solin} leads to several boundary integral
equations defined on the boundary of the scatterer. Four of standard
boundary integral equations on $\Gamma \times (t>0)$ are given as follows:
\begin{align}
\label{out1}
u^{\mathrm{inc}}(x,t)-S^1q(x,t) =0&  \\
\label{out2}
\dot{u}^{\mathrm{inc}}(x,t)-\dot{S}^1q(x,t) =0&  \\
\label{out2.5}
\frac{\partial u^{\mathrm{inc}}}{\partial n}(x,t)-(D^{T1})^-q(x,t)=0&\\
\label{out3}
\frac{\partial u^{\mathrm{inc}}}{\partial n}(x,t)+\frac{1}{c_1}\dot{u}^{\mathrm{inc}}(x,t)-(D^{T1})^-q(x,t)-&\frac{1}{c_1}\dot{S}^1q(x,t) =0
\end{align}
where $\dot{\phantom{u}}$ stands for the time derivative.  Equations
of these types have been considered by many authors for various
potentials mainly in 3D.  Indeed, \eqref{out1} is the ordinary
BIE. The time differentiated equation in \eqref{out2} has been
considered in \cite{mansur}. Bamberger and Ha Duong\cite{bamberger}
also discussed a fully variational version of this equation in 3D.
Equation similar to \eqref{out2.5} for the double layer potential in
3D has been utilised by Parot et al.\cite{parot1,parot2} while
\eqref{out3} for the double layer potential in 3D has been considered
in Ergin et al.\cite{ergin} and Chappell et al.\cite{chappell} among
others. The coupling constant in \eqref{out3} seems to be the most
natural choice because \eqref{out3} is then derived as one imposes the
first order approximation of the absorbing boundary
condition\cite{engquist} on $\Gamma$ to the RHS (right hand side) of
\eqref{solin}, thus ``exteriorising'' the interior domain $D_2$. The
discussion in Chappell and Harris\cite{chappell2} also seems to support
this choice. We shall, however, return to this issue later.

\subsection{Stability}
We consider the following Volterra integral equation which is typically a time domain BIE obtained by discretising \eqref{out1}--\eqref{out3} in the spatial direction by using collocation:
\begin{align}
\label{BIE}
f(t)=\int^t_0 K(t-s)v(s) \ ds
\end{align}
where $K$ represents an $N\times N$ matrix and $v$ and $f$ stand for unknown and given $N$-vectors, respectively. Note that $K$ may include terms of the 
form $c \delta(\cdot)$ or its derivatives, where $c$ is a constant and $\delta(\cdot)$ is Dirac's delta functions.
Discretising the unknown function $v(s)$ in \eqref{BIE} using time interpolation functions $\phi_m(s)$ as $\displaystyle v(s)\approx \sum_m \phi_m(s) v_m$, we obtain the following algebraic equation:
\begin{align}
\label{alge}
f(l\Delta t)& = \sum^l_{m=1}\int^{l\Delta t}_0 K(l\Delta t-s)\phi_m(s) \ ds \ v_m\\
\phi_m(s)&=\phi_{\Delta t}(s-m\Delta t)\notag
\end{align}
where $\phi_{\Delta t}(t)$ is a basis function which satisfies $\phi_{\Delta t}(k \Delta t)=\delta_{k0}$ (where $\delta_{ij}$ is the Kronecker delta),
$\Delta t$ is the time increment, $l$ is the number of time steps, respectively. Usually, an algebraic equation in the form of \eqref{alge} is solved in a time marching manner for $v_m \, (m=1,2,\cdots)$ in time domain BIEMs.

Obviously, the stability of the resulting numerical scheme is a concern
in solving BIEs in time domain. To examine this issue, we follow the standard
argument\cite{walker} to put $v_m=\lambda^m v$
in the homogeneous version of \eqref{alge} where $\lambda\in \mathbb{C}$ is a number and $v$ is an element of $\mathbb{C}^N$. This gives
\begin{align}
\label{algeh}
0 = \sum^{l-1}_{m=0}\int^{l\Delta t}_0 K(s)\phi_{\Delta t}(m\Delta t - s) \ ds \lambda^{-m} v.
\end{align}
Suppose $l$ is taken sufficiently large. A complex number $\lambda$ is
said to be an eigenvalue of \eqref{algeh} if there exists a
non-trivial $v$ which satisfies \eqref{algeh}. Then our definition
of the stability is the following: the scheme is
stable if all the eigenvalues of \eqref{algeh} satisfy
$|\lambda|\le 1$. The scheme is unstable if there exists an eigenvalue
of \eqref{algeh} s.t. $|\lambda|> 1$ holds.  Eigenvalue problems of this
type in 3D have been considered by many authors
after converting them into equivalent linear eigenvalue problems for the
companion matrices (e.g., \cite{walker,wang,parot2,jang2,vantwout,bai}). As a matter of fact, there is no ambiguity in the choice
of a sufficiently large $l$ in 3D if the scatterer is bounded because the fundamental solution has a ``tail'' of a finite length. In 2D problems, however, this is not the case since the tail of the fundamental solution has an infinite length, as one sees in \eqref{eq:fs}. In addition, the time decay of the fundamental solution is slow, thus making it difficult to set an appropriate truncation number $l$ in \eqref{algeh}.

To proceed further, we put $\lambda=e^{-i \Omega \Delta t}$ where
$\Omega$ is a complex number. The stability criterion is now rewritten
as follows:  $\mathrm{Im}\, \Omega\le 0$ ($\mathrm{Im}\, \Omega>0$) implies
stability (instability).  Also, suppose that $\phi_{\Delta t}(s)=0$ for
$s < -\Delta t$; a condition satisfied by many choices of the basis
function including a piecewise linear one.  We then let $l$ tend to infinity
in \eqref{algeh} to have
\begin{align}
\label{eq:dft}
0 = \sum^{\infty}_{m=-\infty}\int^{\infty}_0 K(s)\psi_{\Delta t}(s-m\Delta t) e^{i m  \Delta t \Omega}\, ds \, v,
\end{align}
where $\psi_{\Delta t}(s)=\phi_{\Delta t}(-s)$, which is nothing other than
the discretised Fourier transform of $K$. Obviously, this expression approximates the Fourier transform $\hat K$ of $K$  precisely in lower frequencies, but
just roughly in higher frequencies. This suggests a connection between
the stability of the time domain BIEM and the eigenvalues of the frequency
domain BIEM; an observation made by many authors (e.g., \cite{ergin,chappell}).

We now write $K$ in terms of $\hat K$ as
\begin{align*}
  K(s)= \frac 1{2\pi} \int_{-\infty}^{\infty} \hat K(\omega) e^{-i\omega s} \,d\omega.
\end{align*}
Using the Poisson summation formula, we rewrite \eqref{eq:dft} into
\begin{align}
\label{stab}
0 = \sum^\infty_{m=-\infty} \frac{1}{\Delta t} \hat{K}(\Omega_m)\hat{\phi}_{\Delta t}(\Omega_m) \ v,\quad \Omega_m=\Omega- \frac{2m\pi}{\Delta t}
\end{align}
where $\hat{\phi}_{\Delta t}$ is the Fourier transform of
$\phi_{\Delta t}$ which is given by
\begin{align*}
  \hat{\phi}_{\Delta t}(\Omega)=\frac 2{\Omega^2 \Delta t}(1-\cos \Omega \Delta t)
\end{align*}
for the particular case of the piecewise linear ${\phi}_{\Delta t}$.
The stability issue of the time domain BIEM is
thus reduced to a non-linear eigenvalue problem of finding
$\Omega \in \mathbb{C}$ with which \eqref{stab} has a non-trivial
solution $v\in \mathbb{C}^N$.  Hence, we call these eigenvalues
$\Omega$'s as the characteristic roots of \eqref{stab}. We note that the expression
on the right hand side of \eqref{stab} is periodic with respect to $\Omega$ with the
period of $2\pi/\Delta t$.

We now consider the limit of $\Delta t \to 0$ in \eqref{stab} in a
somewhat intuitive manner. More rigorous arguments could be made with
particular choices of kernel and basis functions. We first note
\begin{align*}
\int_{-\infty}^{\infty}\frac{\phi_{\Delta t}(s)}{\Delta t}\, ds=1
\end{align*}
if $\phi_{\Delta t}$ can interpolate a constant function exactly. This implies
\begin{align*}
\frac{\phi_{\Delta t}(s)}{\Delta t}\to \delta(s) \qquad \mbox{as }\Delta t \to 0,
\end{align*}
a conclusion which could be checked with particular choices of $\phi_{\Delta t}$. This gives
\begin{align*}
\lim_{\Delta t \to 0}\frac{\hat\phi_{\Delta t}}{\Delta t}\to 1
\end{align*}
Hence, we expect to have
\begin{align}
\label{freq}
\sum^\infty_{m=-\infty} \frac{1}{\Delta t} \hat{K}(\Omega_m)\hat{\phi}_{\Delta t}(\Omega_m) \to \hat{K}(\Omega) \quad \mbox{as }\Delta t \to 0,
\end{align}
if $\hat K(\Omega) \to 0$ as $|\Omega| \to \infty$, which 
is the case in 2D. From this result, we
expect that the characteristic roots of \eqref{stab} are obtained as
perturbations of the eigenvalues of the corresponding frequency domain
BIEs. We note that a similar observation has been made in Chappell et
al.\cite{chappell,chappell2} qualitatively.

It is well-known that the eigenvalues of the frequency domain BIE can
be classified into true and fictitious
eigenvalues\cite{misawa1,misawa2}. In the exterior problems, the true
eigenvalues are with negative imaginary parts, while the behaviour of
the fictitious eigenvalues depend on the particular choice of integral
equations. In \eqref{out1}--\eqref{out2.5} the fictitious eigenvalues
of the corresponding frequency domain BIEs are real valued, while
those of \eqref{out3} are with negative imaginary parts. It is
therefore natural to expect that equations
\eqref{out1}--\eqref{out2.5} are more prone to instability than
\eqref{out3}. However, \eqref{out1}--\eqref{out2.5} may still turn out
to be stable after discretisation depending on the choice of the time
basis function because real eigenvalues of the frequency domain BIE
may move to lower complex plane after the time discretisation. Also,
\eqref{out3} may turn out to be unstable if the perturbation of the
eigenvalues is very large.

\subsection{Simplified stability analysis for circular domains}
\label{fourier}
One may use methods based on contour integrals such as the
Sakurai-Sugiura Method (SSM) in the solution of non-linear eigenvalue
problem in \eqref{stab} for a general boundary $\Gamma$.
Indeed, one may replace 
\begin{align*}
\hat G^{\nu}(x)=\frac{i}4 H_0^{(1)}(\Omega |x|/c^\nu)
\end{align*}
in the Fourier transformed versions of BIEs in \eqref{out1}--\eqref{out3} by
\begin{align}
\label{eq:ftg}
\frac{i}4 \sum^\infty_{m=-\infty} H_0^{(1)}(\Omega_m|x|/c^\nu)\frac{\hat{\phi}_{\Delta t}(\Omega_m)}{\Delta t}
\end{align}
to this end, where $H_0^{(1)}$ is the Hankel function of the 1st kind.  In the
present paper, however, we shall pay attention to a simple special case
in which $\Gamma$ is a unit circle. Also, we restrict out attention to
the piecewise linear time basis functions for the purpose of
simplicity.

We consider \eqref{stab} for a unit circle $\Gamma$ without spatial
discretisation (the effect of the spatial discretisation will be
discussed later). In this case we can simplify the non-linear
eigenvalue problem in \eqref{stab} using the Fourier series
with respect to the angular variable. Indeed, we use the well-known
Graf addition theorem\cite{martin} to have
\begin{align}
\label{graf}
  \hat G(x-y)=\frac{i}4 \sum_{n=-\infty}^{\infty}H_n^{(1)}(k|x|)J_n(k|y|) e^{i n (\Theta - \theta)}, &\\ x=|x|(\cos \Theta, \sin \Theta),\quad y=|y|(\cos \theta, \sin \theta)&\notag
\end{align}
when $|x|>|y|$ holds, where $J_n$ is the Bessel function, $k=\Omega/c$
and $\Theta$ and $\theta$ are the azimuth angles of $x$ and $y$. The
role of $H_n^{(1)}$ and $J_n$ in \eqref{graf} is interchanged when
$|x|<|y|$. In \eqref{graf} we have suppressed the superfix $\nu$ for the
domain in order to simplify the notation.  From \eqref{stab} and
\eqref{graf}, we see that the characteristic roots of the time
discretised boundary integral equations corresponding to
\eqref{out1}--\eqref{out2.5} are obtained as zeros of the expressions
in the following list:
 \begin{align}
\label{stabs}
  &  S\leftrightarrow \sum_{m}H_n^{(1)}(\Omega_m/c)J_n(\Omega_m/c)\hat{\phi}_{\Delta t}(\Omega_m) \\
\label{stabds}
  &  \dot{S}\leftrightarrow f_1(\Omega; n, c)=-\sum_{m}i\Omega_m H_n^{(1)}(\Omega_m/c)J_n(\Omega_m/c)\hat{\phi}_{\Delta t}(\Omega_m) \\
\label{stabdt}
  &  D^{T-}\,(,D^+)\leftrightarrow f_2(\Omega; n, c)=\sum_{m} \Omega_m/cH_n^{(1)}(\Omega_m/c)J'_n(\Omega_m/c)\hat{\phi}_{\Delta t}(\Omega_m)
\end{align}
where $n$ is an integer between 0 and a large number $n_{\mathrm{max}}$.
The characteristic equation for \eqref{out3} is
obtained from \eqref{stabds} and \eqref{stabdt} as
\begin{align}
\label{stabbm}
  \sum_{m} H_n^{(1)}(\Omega_m/c) \Omega_m/c(J'_n(\Omega_m/c)-i J_n(\Omega_m/c))\,\hat{\phi}_{\Delta t}(\Omega_m).
\end{align}
Note that the series on the right hand sides of eqs. \eqref{stabs}--\eqref{stabbm} are absolutely convergent. 

\subsection{Numerical experiments}
\label{num_exp_dirichlet}
We now carry out numerical experiments to see if the 
stability analysis given in the previous section can predict the behaviour
of the time domain BIEM correctly. 

To this end, we consider the problem defined in \ref{sec:edp} where the boundary $\Gamma$ is the unit circle. The material constants are $s_1=\rho_1=c_1=1$.
The incident wave is a plane wave given by:
\begin{align}
\label{quadinc}
u^{\mathrm{inc}}=
\begin{cases}
0 & (c_1t-x_1-t_0 \le 0) \\
\frac{(c_1t-x_1-t_0)^2}{2} & (c_1t-x_1-t_0 > 0)
\end{cases}
\end{align}
where $t_0=1+2\Delta t$.  We use piecewise constant boundary elements,
piecewise linear temporal elements and the collocation method to
discretise the BIEs in \eqref{out1}--\eqref{out3}.  All the required
integrals are computed exactly.  The boundary is discretised into 100
elements, the time increment is set as $\Delta t=\frac{2\pi}{100}$ and
the number of time steps is 1000.  Also, the characteristic roots of
\eqref{stab} are calculated with \eqref{stabs}--\eqref{stabbm} and
SSM.

\begin{figure}[htb]
  \begin{center}
\includegraphics[width=\linewidth]{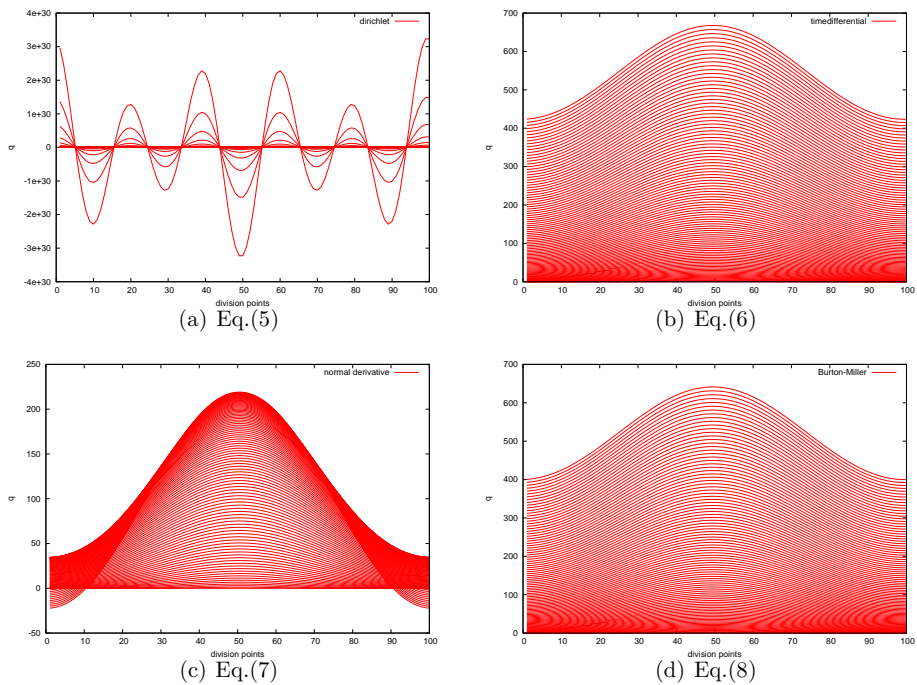}
    \caption{$q$ obtained with various integral equations vs point number. The point number and the azimuth angle $\theta$ are related by $\theta=-2\pi/100\times{}$ point number} 
  \label{fig:1}
  \end{center}
         \end{figure}

\begin{figure}[htb]
\begin{center}
           \includegraphics[width=0.6\linewidth]{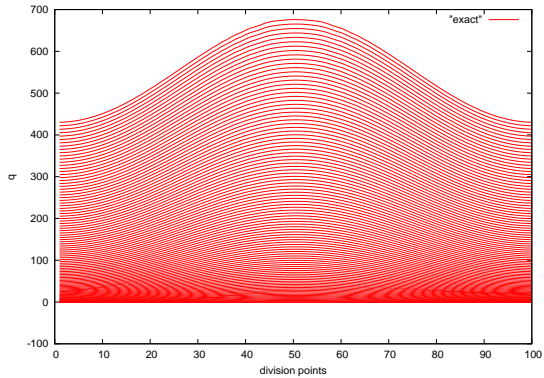}
 \end{center}
  \caption{``Exact'' solution of $q$}
  \label{exact_q}
\end{figure}

Figs.\ref{fig:1}(a)--\ref{fig:1}(d) show the results obtained with
\eqref{out1}--\eqref{out3}, respectively. We plot $q$ for every 10
time steps in these figures (this applies to all subsequent time
domain results).  Also, Fig.\ref{exact_q} gives the ``exact'' solution
obtained numerically with the frequency domain exact solution and
FFT. We see that the standard BIE in \eqref{out1} is unstable, and the
time derivative BIE in \eqref{out2} and the time domain BM BIE in
\eqref{out3} are stable. The normal derivative BIE in \eqref{out2.5}
does not show divergence, but deviates considerably from the ``exact''
solution. The BM result is not as bad as the normal derivative one,
but is seen to drift from the exact solution by a time dependent
constant. The accuracy of the time derivative BIE appears to be
satisfactory.

We next check the behaviour of the characteristic roots of these time
domain BIEs using SSM. We set the range for computing eigenvalues to
be $0 \le \mathrm{Re} \ \Omega \le 50=\frac{\pi}{\Delta t}$ and
$-2 \le \mathrm{Im} \ \Omega \le 2$ considering the periodicity of
\eqref{stab} and the fact that the characteristic roots are located
symmetrically with respect to the imaginary axis, which can be easily
shown using the explicit forms of \eqref{stab}. Note that the upper
limit of $\Omega$ is consistent with the Nyquist frequency associated
with $\Delta t$.  Also, we took $n_{\mathrm{max}}$ to be 60 considering the
number of boundary subdivisions and the spatial Nyquist ``frequency''.
In the computation, we redefine the Hankel functions so that they have
branch cuts on the negative imaginary axis rather than on the negative
real axis. This guarantees that the expression in \eqref{stab} is
analytical when $0<\mathrm{Re}\,\Omega<\frac{\pi}{\Delta t}$ holds.

\begin{figure}[htb]
  \begin{center}
\includegraphics[width=\linewidth]{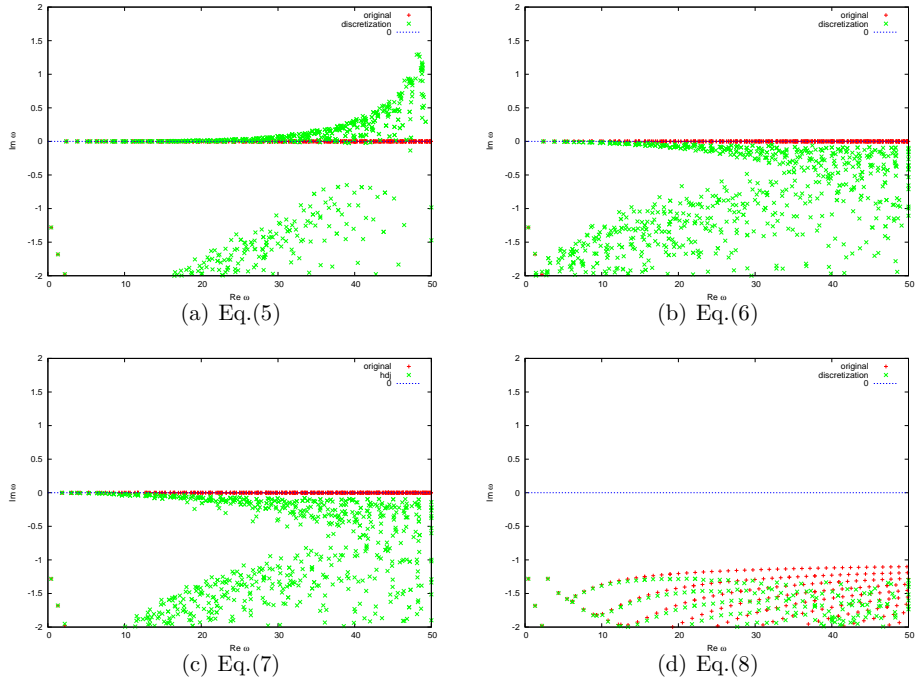}
  \end{center}
  \caption{Characteristic roots of various integral equations}
  \label{fig:eig}
  \end{figure}

Figs.\ref{fig:eig}(a)--Fig.\ref{fig:eig}(d) show the characteristic roots
of the BIEs in \eqref{out1}--\eqref{out3}, respectively. We plot the
eigenvalues of \eqref{stab} for various BIEs (i.e., zeros of the
expressions in \eqref{stabs}--\eqref{stabbm}) in green and the
eigenvalues of the frequency domain BIEs given by $\hat K(\Omega) v =0$ (i.e.,
zeros of the products of Hankel and Bessel functions obtained by
setting $m=0$ in \eqref{stabs}--\eqref{stabbm}) in red.
  Note that the red symbols near the imaginary axis in
  Figs.\ref{fig:eig}(a)--Fig.\ref{fig:eig}(d) are the true eigenvalues for
  the exterior Dirichlet problem, while those on the real axis are
  fictitious ones related to the interior Dirichlet problems in
  Figs.\ref{fig:eig}(a) and \ref{fig:eig}(b) and to the
  Neumann problem in Fig.\ref{fig:eig}(c). The fictitious eigenvalues
  for \eqref{out3} are those associated with interior impedance
  boundary value problems. We note that all BIEs in
  \eqref{out1}--\eqref{out3} have characteristic roots close to the
  true or fictitious eigenvalues, but other characteristic roots are
  scattered and quite far from any of eigenvalues of the corresponding
  frequency domain BIEs, except in the BM equation in \eqref{out3}.
      
From these figures, we see that the BIE \eqref{out1} has the
characteristic roots with positive imaginary parts, but this is not
the case with other BIEs.  These results are consistent with the
corresponding time domain results.  Also, the inaccuracy of
\eqref{out2.5} is considered to be related to the fact that $\Omega=0$
is an eigenvalue of \eqref{stabdt} with $n=0$. We remark that a
similar case has been reported in Parot et al.\cite{parot1} where this
phenomenon has been called a ``pneumatic mode''.
As a matter of fact $\Omega=0$ is a zero of both \eqref{stabds} (for
all $n$) and \eqref{stabbm} (for $n=0$) as well. An adverse effect of
this eigenvalue on \eqref{out3} is visible in the constant shift of
the solution in Fig.\ref{fig:1}(d), although not as evidently as in
Fig.\ref{fig:1}(c). 

\begin{figure}[htb]
  \begin{center}
    \includegraphics[width=\linewidth]{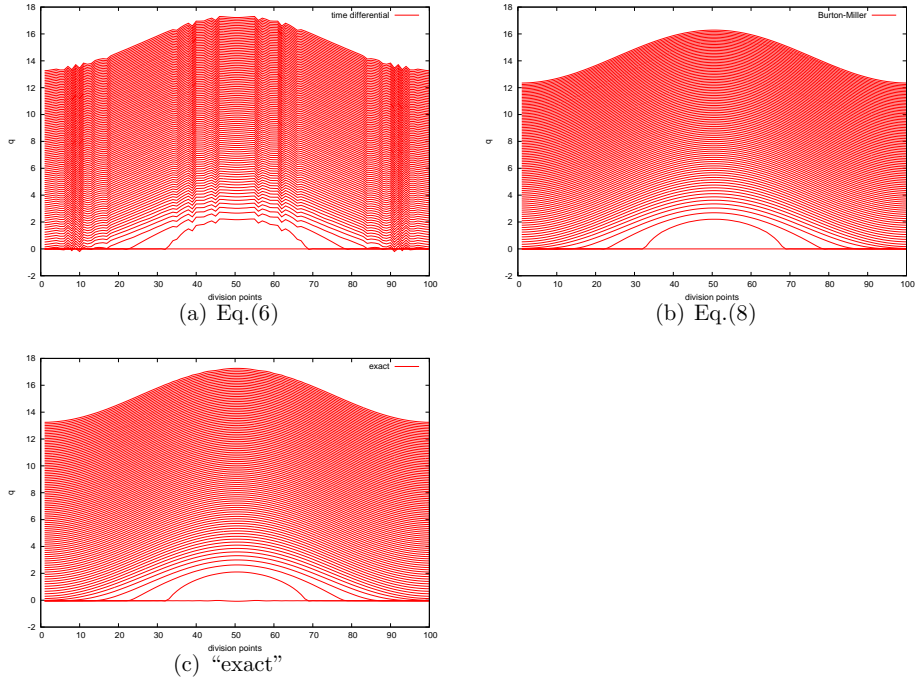}
    \caption{$q$ obtained with various integral equations}
    \label{fig:linear}
  \end{center}
\end{figure}

To examine the effect of this zero eigenvalue on the numerical
solution of \eqref{out2}, we consider another incident wave given by
\begin{align}
\label{lininc}
u^{\mathrm{inc}}=
\begin{cases}
0 & (c_1t-x_1-t_0 \le 0) \\
\frac{(c_1t-x_1-t_0)^2}{c_1t-x_1-t_0 + 4 \Delta t} & (c_1t-x_1-t_0 > 0),
\end{cases}
\end{align}
which is a smoothed linear function, where $t_0=1+2\Delta t$. Setting
other parameters the same as in the previous example, we solve
\eqref{out2} to compute the time history of $q$ as plotted in
Fig.\ref{fig:linear}(a). Comparing this result with the ``exact'' solution
given in Fig.\ref{fig:linear}(c) one sees that an error having a zig-zag
pattern is superimposed on the solution of \eqref{out2}. This is in
contrast to the BM solution given in Fig.\ref{fig:linear}(b) which is
smooth. This result can be explained as follows.  With \eqref{out2},
the (spatially) high frequency error incurred initially by the
mismatch of the wavefront and mesh remains undamped after a long time
because of the existence of a zero characteristic root with high $n$
eigenfunctions. Since this eigenvalue is zero, this error does not
propagate, decay or amplify. In other words, it persists. This type of
error is included also in Fig.\ref{fig:1}(b), although its magnitude is
too small to be visible. From these numerical experiments, we conclude
that none of the integral equations in \eqref{out1}--\eqref{out3} are
satisfactory!

\begin{figure}[htb]
  \begin{center}
\includegraphics[width=\linewidth]{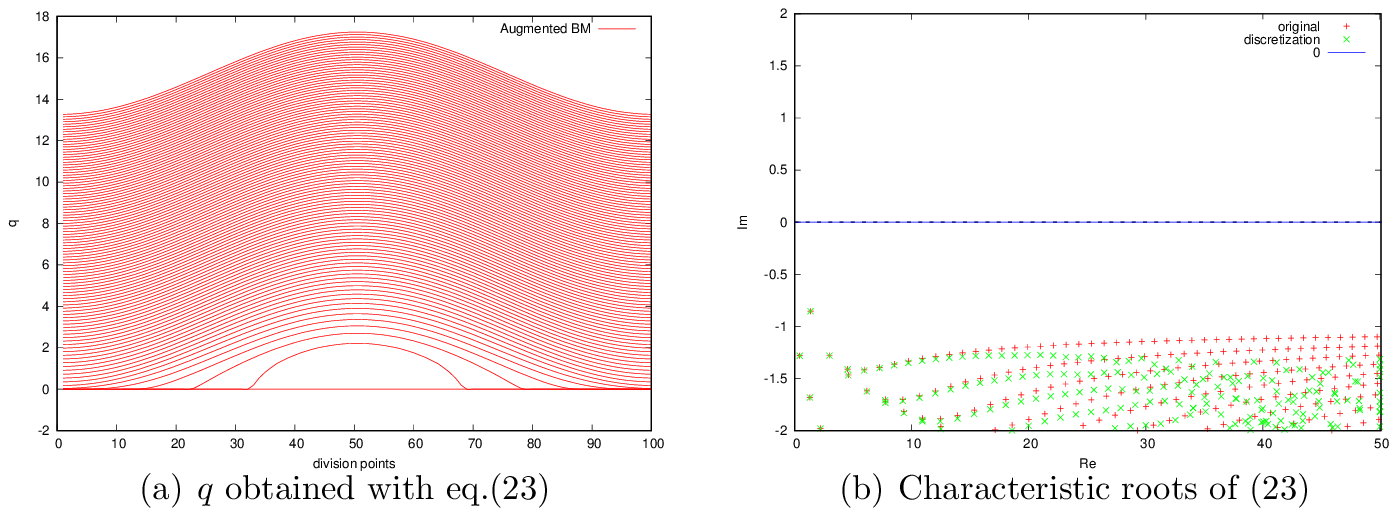}
\caption{Augmented BM integral equation}
\label{fig:modBM}
  \end{center}
\end{figure}

A possible remedy for all these problems is to use an integral equation
given by
\begin{align}
\label{out4}
\frac{\partial u^{\mathrm{inc}}}{\partial n}(x,t)&+\frac{1}{c_1}\dot{u}^{\mathrm{inc}}(x,t)+ \alpha {u}^{\mathrm{inc}}(x,t)\notag\\
&-(D^{T1})^-q(x,t)-\frac{1}{c_1}\dot{S}^1q(x,t) -\alpha S^1q(x,t) =0,
\end{align}
which is the time domain counterpart of the BM equation with a complex
(not pure imaginary) coupling constant, where $\alpha$ is a (real)
number. It is easy to see that $\Omega=0$ is not a characteristic
root of this equation. The numerical solution $q$ obtained with
\eqref{out4} and the incident wave in \eqref{lininc} is given in Fig.\ref{fig:modBM}(a) and its characteristic roots are shown in
Fig.\ref{fig:modBM}(b), where we set $\alpha=1$.
The high accuracy and stability of this formulation is evident from these
figures.

Epstein et al.\cite{epstein} investigated time domain integral
equations whose solutions exhibit correct long-time behaviours.  From
this viewpoint, \eqref{out4} seems to be a better choice than other
stable choices in \eqref{out2}--\eqref{out3}, although the
characteristic root of \eqref{out4} whose imaginary part has the
minimum magnitude is close to one of fictitious eigenvalues (the one
whose imaginary part is approximately equal to -0.8) rather than a
true one.

\section{``Stable potentials''}
\label{stable_potentials}
\begin{figure}[htb]
 \begin{center}
\includegraphics[width=\linewidth]{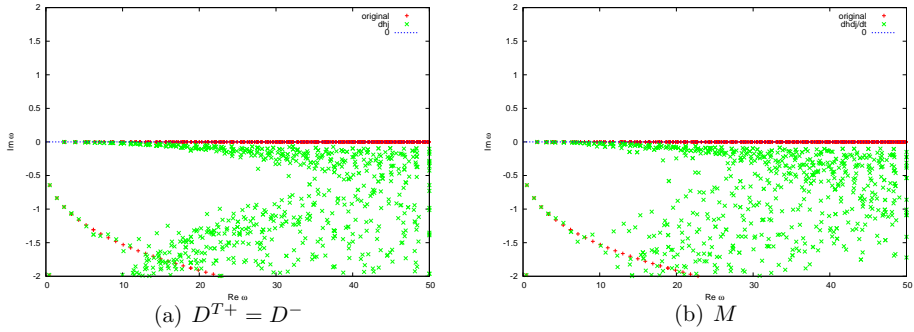}
\caption{Characteristic roots of various integral operators}
\label{fig:6}
 \end{center}
\end{figure}

Motivated by the results in the previous section, we examine the
stability of integral equations on the unit circle derived from
potentials which may appear in BIEs. These potentials include the
single layer $S\, (=S^+=S^-)$, traces of the normal derivatives of $S$
denoted by $D^{T\pm}$, the time derivative of $S$ denoted by
$\dot S\, (=\dot S^+=\dot S^-)$ and the traces of the double layer
$D^{\pm}$. Although we are also interested in the normal derivative of
$D$ denoted by $N$, it turned out that the simplified approach
presented in \ref{fourier} using the Fourier series expansion is
not very easy to apply to $N$ with piecewise linear time basis
functions because the series similar to \eqref{stabs}--\eqref{stabdt}
for $N$ does not converge absolutely (A similar observation applies to
$\dot D^\pm$ as well).  Using a smoother time basis function could be
a solution.  As we shall see later, however, the time integrated
normal derivative of the double layer potential defined by
\begin{align*}
  M u = \int_0^t\int_{\Gamma}\frac{\partial}{\partial n_x}\frac{\partial}{\partial n_y}
  \log \frac{c (t-s) + \sqrt{c^2 (t-s)^2 - |x-y|^2}_+}{|x-y|} u(y,s) dS_y ds\\
x \in \mathbb{R}^2\setminus \Gamma,\quad t>0
\end{align*}  
is more useful than $N$ as far as the stability is concerned. We
therefore carry out the stability analysis in \ref{fourier} with
$M \, (=M^+=M^-)$ instead of $N$.  Since the results for $S$, $\dot S$ and
$D^{T-}=D^+$ have already been given in
Figs.\ref{fig:eig}(a)--\ref{fig:eig}(c), we present those for $D^{T+}=D^-$
and $M$ in Figs.\ref{fig:6}(a)--\ref{fig:6}(b) using the same time
increment as before (i.e., $\Delta t=2\pi/100$).
They are zeros of the following expressions.
 \begin{align}
\label{Dplus}
  &  D^{-},\,D^{T+}\leftrightarrow f_3(\Omega; n, c)=\sum_{m}\Omega_m/c  H_n^{(1)'}(\Omega_m/c)J_n(\Omega_m/c)\hat{\phi}_{\Delta t}(\Omega_m)\\
\label{M_}
  &  M\leftrightarrow f_4(\Omega; n, c)=\sum_{m} i\Omega_m/c^2 H_n^{(1)'}(\Omega_m/c)J'_n(\Omega_m/c) \hat{\phi}_{\Delta t}(\Omega_m)
\end{align}

From these results, we confirm that the equations obtained by
discretising the following integral equations are numerically stable
with piecewise linear time basis functions: (a) the time derivative of
the single layer potential (b) the interior and exterior traces of the
normal derivative of the single layer potential (c) the interior and
exterior traces of the double layer potential (d) the time integrated
normal derivative of the double layer potential. Note, however, that
we have no claim of stability of these potentials except in the cases
tested here.

In the next section, we combine these ``stable potentials'' to
obtain numerically stable formulations in transmission problems.

\section{Transmission problems}

We are now interested in finding $u$ which satisfies \eqref{wave},
\begin{align*}
\Delta u -\frac{1}{c_2^2}\frac{\partial^2 u}{{\partial t}^2} = 0 \hspace{0.2truecm} \mathrm{in} \  D_2 \times (t>0),
\end{align*}
the transmission boundary conditions on $\Gamma$:
\begin{align*}
u^+=u^-(=u),\quad
s_1\frac{\partial u^+}{\partial n} = s_2\frac{\partial u^-}{\partial n} \ (=q) \quad
\mbox{ on }\Gamma \times (t>0)
\end{align*}
and the homogeneous initial conditions
\begin{align*}
  u|_{t=0}=\dot u|_{t=0}=0\qquad \mbox{in }D_2
\end{align*}
in addition to the homogeneous initial and radiation conditions for
$u^{\mathrm{sca}}$ in \eqref{eq:init}, where $c_2$ is the wave speed
in $D_2$ given by $c_2=\sqrt{\frac{s_2}{\rho_2}}$ and ($s_2$,
$\rho_2$) are the shear modulus and density in $D_2$,
respectively. The superscript $+(-)$ stands for the trace to $\Gamma$
from $D_1$ ($D_2$), respectively.

\subsection{Boundary integral equations}
There exist various possibilities of integral equation for
transmission problems on $\Gamma$, of which we consider the following
four \cite{misawa1,misawa2}:

PMCHWT
\begin{align}
\label{PMold}
  \left(
    \begin{array}{cc}
      -(D^{1}+D^{2}) & \frac{1}{s_1}S^{1}+\frac{1}{s_2}S^{2} \\
      -(s_1N^{1}+s_2N^{2}) & D^{T1}+D^{T2}  
    \end{array}
  \right)
\left(
\begin{array}{c}
u \\
q
\end{array}
\right)
=\left(
    \begin{array}{c}
    u^{\mathrm{inc}}  \\
    s_1\frac{\partial u^{\mathrm{inc}}}{\partial n}  
    \end{array}
  \right)
\end{align}

M\"uller
\begin{align}
\label{Mold}
    \left(
    \begin{array}{cc}
      \frac{s_1+s_2}{2}-(s_1D^1-s_2D^2) & S^1-S^2 \\
      -(N^1-N^2) & \frac{s_1+s_2}{2s_1s_2}+\frac{1}{s_1}D^{T1}-\frac{1}{s_2}D^{T2}  
    \end{array}
  \right)
\left(
\begin{array}{c}
u \\
q
\end{array}
\right)
=
\left(
\begin{array}{c}
s_1u^{\mathrm{inc}}\\
\frac{\partial u^{\mathrm{inc}}}{\partial n}
\end{array}
\right)
\end{align}

Burton-Miller
\begin{align}
\label{BMold}
    \left(
    \begin{array}{cc}
       \frac{1}{2c_1}\frac{\partial}{\partial t}-(N^1+\frac{1}{c_1}\dot{D^1}) & \frac{1}{2s_1}+\frac{1}{s_1}D^{T1}+\frac{1}{c_1s_1}\dot{S^1} \\
          -\frac{1}{2}-D^2 & \frac{1}{s_2}S^2
          \end{array}
  \right)
\left(
\begin{array}{c}
u \\
q
\end{array}
\right)
=
\left(
\begin{array}{c}
\frac{\partial u^{\mathrm{inc}}}{\partial n}+\frac{1}{c_1}\frac{\partial u^{\mathrm{inc}}}{\partial t} \\
0
\end{array}
\right)
\end{align}

standard
\begin{align}
\label{sold}
  \left(
    \begin{array}{cc}
      \frac12-D^1 & \frac{1}{s_1}S^1 \\
           \frac12+D^2 & -\frac{1}{s_2}S^2
     \end{array}
  \right)
\left(
\begin{array}{c}
u \\
q
\end{array}
\right)
=\left(
    \begin{array}{c}
    u^{\mathrm{inc}}  \\
0
    \end{array}
  \right)
\end{align}
In these equations we write $D^{\nu}$ for $(D^{\nu +}+D^{\nu
  -})/2$ etc.\ in order to show the non-integral terms explicitly at the
cost of an abuse of notation. We note that there exist several versions
of M\"uller's formulations for the wave equation. We here use the one
in which the singularities of single layer potentials cancel.

\subsection{Stable formulations}
\label{stabon}

The boundary integral equations shown in the previous section
can be rewritten easily in terms of ``stable potentials''
presented in section \ref{stable_potentials} with the help
of time differentiation and integration by parts. Here are the results:

modified PMCHWT
\begin{align}
\label{PMnew}
  \left(
    \begin{array}{cc}
      -(D^1+D^2) & \frac{1}{s_1}\dot S^1+\frac{1}{s_2}\dot S^2 \\
      -(s_1M^1+s_2M^2) & D^{T1}+D^{T2}  
    \end{array}
  \right)
\left(
\begin{array}{c}
\dot u \\
q
\end{array}
\right)
=\left(
\begin{array}{c}
\dot u^{\mathrm{inc}} \\
s_1\frac{\partial u^{\mathrm{inc}}}{\partial n}
\end{array}
\right)
\end{align}

modified M\"uller
\begin{align}
\label{Mnew}
    \left(
    \begin{array}{cc}
      \frac{s_1+s_2}{2}-(s_1D^1-s_2D^2) & \dot S^1-\dot S^2 \\
      -(M^1-M^2) & \frac{s_1+s_2}{2s_1s_2}+\frac{1}{s_1}D^{T1}-\frac{1}{s_2}D^{T2}  
    \end{array}
  \right)
\left(
\begin{array}{c}
\dot u \\
q
\end{array}
\right)
=
\left(
\begin{array}{c}
s_1 \dot u^{\mathrm{inc}}\\
\frac{\partial u^{\mathrm{inc}}}{\partial n}
\end{array}
\right)
\end{align}

modified Burton-Miller
\begin{align}
\label{BMnew}
    \left(
    \begin{array}{cc}
       \frac{1}{2c_1}-(M^1+\frac{1}{c_1}{D^1}) & \frac{1}{2s_1}+\frac{1}{s_1}D^{T1}+\frac{1}{c_1s_1}\dot{S^1} \\
          -\frac{1}{2}-D^2 & \frac{1}{s_2}\dot S^2
          \end{array}
  \right)
\left(
\begin{array}{c}
\dot u \\
q
\end{array}
\right)
=
\left(
\begin{array}{c}
\frac{\partial u^{\mathrm{inc}}}{\partial n}+\frac{1}{c_1}\frac{\partial u^{\mathrm{inc}}}{\partial t} \\
0
\end{array}
\right)
\end{align}

modified standard
\begin{align}
\label{snew}
  \left(
    \begin{array}{cc}
      \frac12-D^1 & \frac{1}{s_1}\dot S^1 \\
           \frac12+D^2 & -\frac{1}{s_2}\dot S^2
     \end{array}
  \right)
\left(
\begin{array}{c}
\dot u \\
q
\end{array}
\right)
=\left(
    \begin{array}{c}
   \dot u^{\mathrm{inc}}  \\
0
    \end{array}
  \right)
\end{align}
where $M^\nu$ is the time integral of $N^\nu$.

The PMCHWT, M\"uller and BM formulations are known not to
have real fictitious eigenfrequencies in the frequency domain, while
the standard formulation does have real fictitious
eigenfrequencies\cite{misawa1,misawa2}. It is therefore expected that
the standard formulation is more prone to instability.

We remark that the time differentiated standard integral equation in
the modified standard equation \eqref{snew} has appeared in the paper by
Panagiotopoulos and Manolis\cite{panag} on elastodynamics in 3D. Also,
the combined use of of $D$, $\dot S$, $M$ and $D^T$ in \eqref{PMnew},
\eqref{Mnew} and \eqref{BMnew} has been proposed by Abboud et
al.\cite{abboud} and Banjai et al.\cite{banjai} in different
contexts in 3D. Their choices of unknowns are different from ours. To the
best of our knowledge, however, these potentials have not been
utilised in forms given in \eqref{PMnew}, \eqref{Mnew} and
\eqref{BMnew} in transmission problems for the wave equation in 2D.

\subsection{Numerical experiments}
\label{ne_transmission}
Setting $s_1=1$, $\rho_1=1$, $s_2=0.2$ and $\rho_2=0.37$, we solve the
time domain BIEs in \eqref{PMold}--\eqref{snew}. The incident wave is
the quadratic one in \eqref{quadinc} and the number of boundary
subdivisions, the number of time steps and the time increments are the
same as in \ref{num_exp_dirichlet}. We use piecewise linear time basis
functions for ($u,q$) in the ordinary formulations and for ($\dot
u,q$) in the modified formulations.

\begin{figure}[htb]
  \begin{center}
\includegraphics[width=\linewidth]{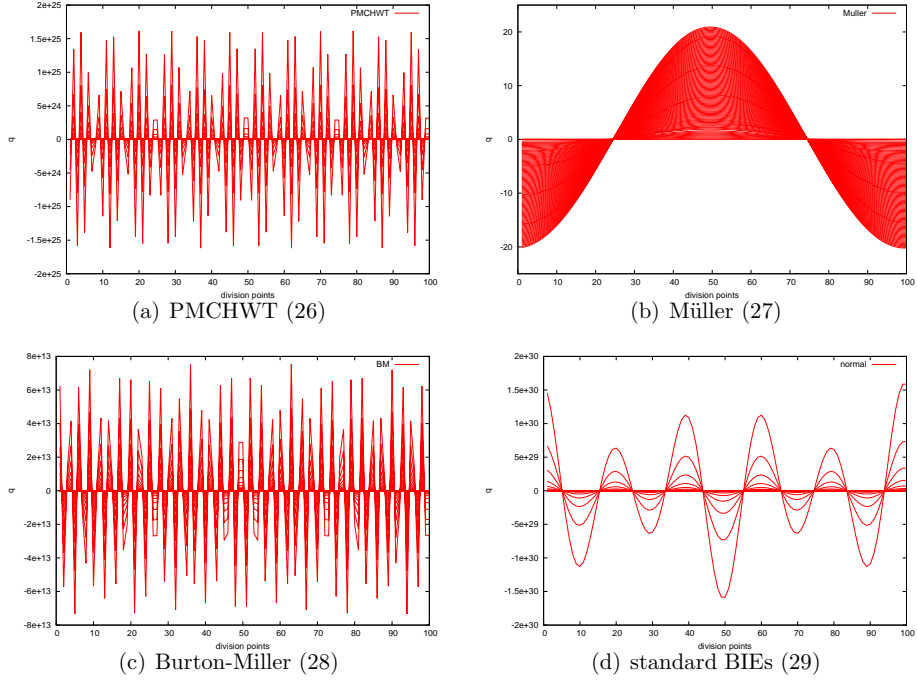}
    \end{center}
\caption{$q$ obtained with ordinary integral equations for transmission problems}
  \label{old}
\end{figure}

\begin{figure}[hbt]
  \begin{center}
\includegraphics[width=\linewidth]{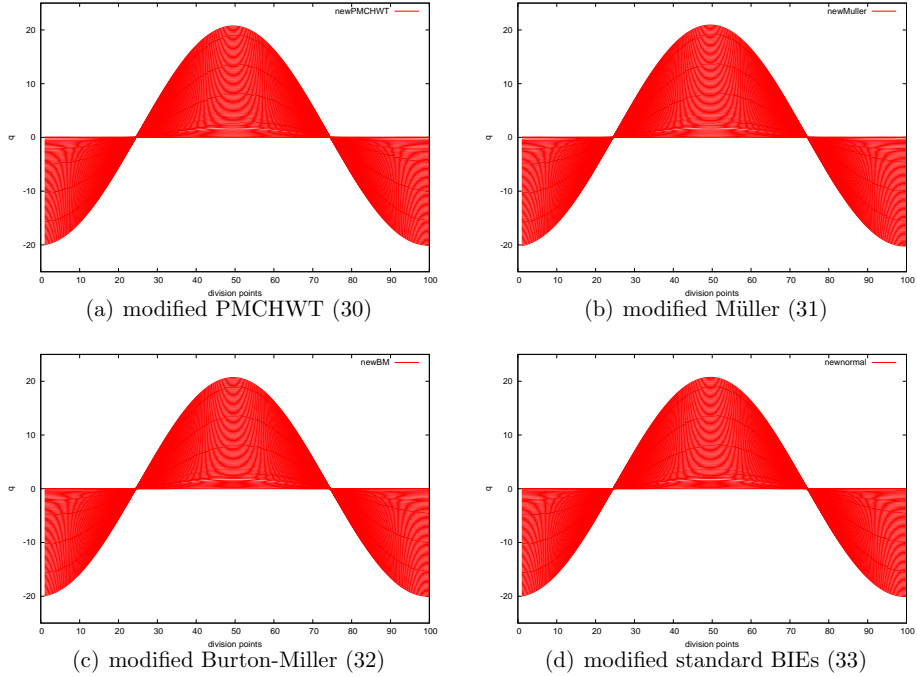}
    \end{center}
\caption{$q$ obtained with modified integral equations for transmission problems}
  \label{new}
\end{figure}

\begin{figure}[htbp]
  \begin{center}
\includegraphics[width=\linewidth]{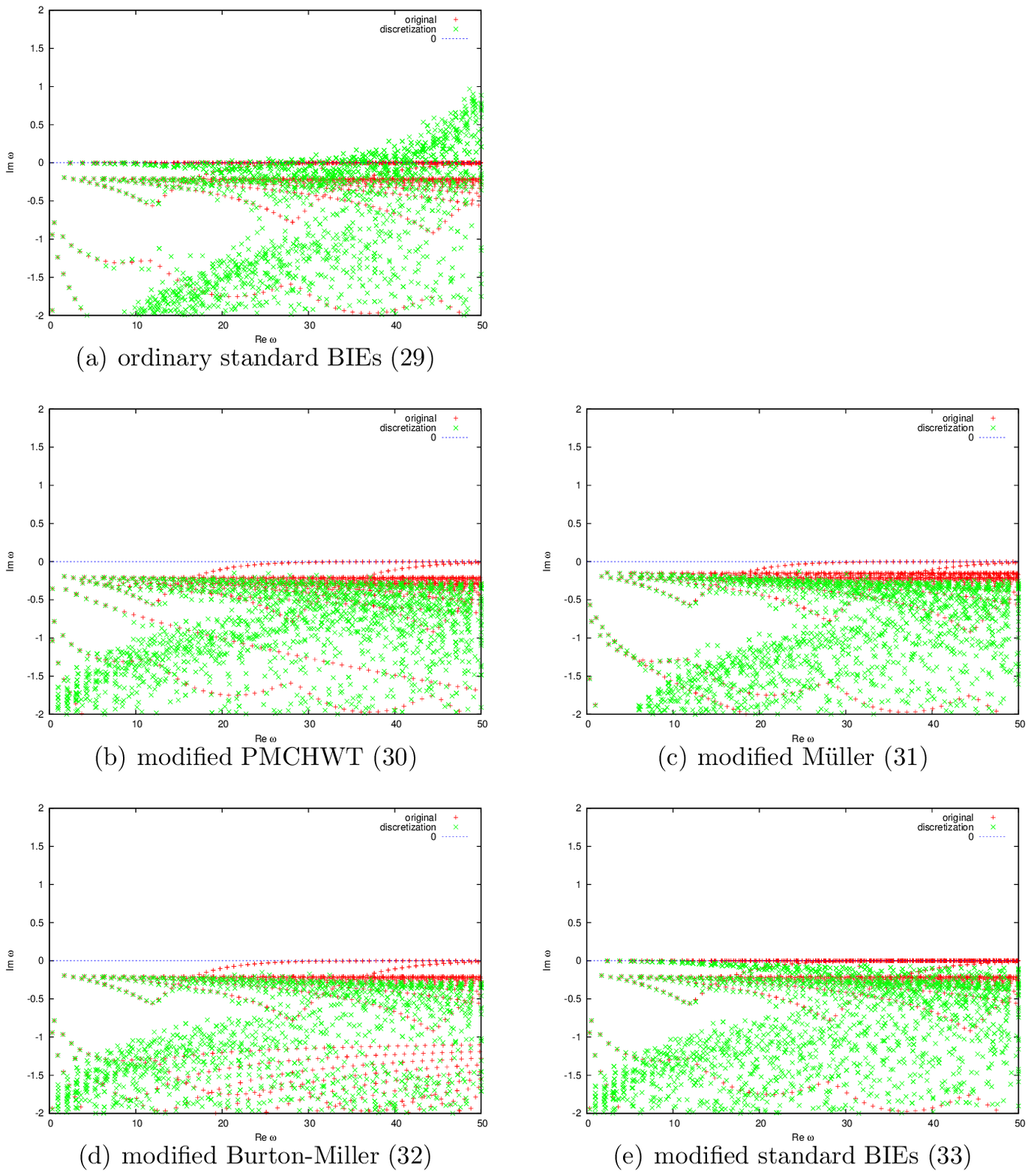}
    \end{center}
\caption{Characteristic roots of various integral equations for transmission problems}
  \label{fig:CR_transmission}
\end{figure}

Fig.\ref{old} and Fig.\ref{new} show the results of $q$ obtained with the ordinary and modified integral equations respectively. We see that the ordinary formulations give unstable results except for the M\"uller formulation, whereas all the modified formulations provide stable results.
Fig.\ref{fig:CR_transmission} shows the distribution of characteristic
roots for those formulations which do not include $N$ or $\dot D$.  In
the modified PMCHWT, for example, they are obtained as the non linear
eigenvalues ($\Omega$'s) of the following matrix for one of
$n=0,\cdots,n_{\mathrm{max}}$. See \eqref{stabds}, \eqref{stabdt}, \eqref{Dplus} and \eqref{M_}: 
\begin{align}
\label{PMnewstab}
 \tilde K&(\Omega;n)_{\mathrm{PMCHWT}}= \notag\\
&\left(
    \begin{array}{cc}
      -(f_3(\Omega;n,c_1)+f_2(\Omega;n,c_2)) & \frac{1}{s_1}f_1(\Omega;n,c_1)+\frac{1}{s_2} f_1(\Omega;n,c_2) \\
      -(s_1f_4(\Omega;n,c_1)+s_2f_4(\Omega;n,c_2)) & f_2(\Omega;n,c_1)+f_3(\Omega;n,c_2)
    \end{array}
  \right)
\end{align}
The distribution of characteristic roots shown in
Fig.\ref{fig:CR_transmission} is seen to be consistent with the time
domain results in Figs.\ref{old} and \ref{new}.  These results also
suggest that the use of standard BIEs may not be recommended even
after the modification since this formulation has many characteristic
roots near the real axis.

Among other three formulations the M\"uller formulation appears to be
better in terms of stability since it gives stable results even
without modification, as one sees in Figs.\ref{old} and
\ref{new}. Another reason to prefer M\"uller is the behaviour of the
characteristic equations \eqref{PMnewstab}, etc., near $\Omega=0$. As a matter
of fact, we can show that $\Omega=0$ is a characteristic root of BM
for $n=0$, but other two do not suffer from this problem.  However,
the characteristic equation \eqref{PMnewstab} for the modified PMCHWT
has the following asymptotic behaviour near $\Omega=0$:
\begin{align}
\label{PMnewstaba}
 \tilde K(\Omega;n)_{\mathrm{PMCHWT}}= \left(
    \begin{array}{cc}
      o(\Omega) & O(\Omega) \\
      O(1/\Omega) & o(\Omega)
    \end{array}
  \right) \qquad n\ge 1
\end{align}
while $\tilde K(\Omega;0)_{\mathrm{PMCHWT}}= O(1)$. This means that the
vector $(0,1)^T$ behaves asymptotically like an eigenvector of
\eqref{PMnewstaba} for $\Omega\approx 0$. This suggests that an arbitrary error
of $q$ having a zero spatial mean may
persist in the solution of the modified PMCHWT. For the
modified M\"uller equation, however, there is no problem of this kind
since we have
\begin{align*}
\tilde K(\Omega;n)_{\text{M\"uller}}= \left(
    \begin{array}{cc}
      O(1) & o(1) \\
      o(1) & O(1)
    \end{array}
  \right).
\end{align*}
We thus conclude that the modified M\"uller equation is a better
choice than the other two in the cases tested.

\begin{figure}[htb]
  \begin{center}
\includegraphics[width=\linewidth]{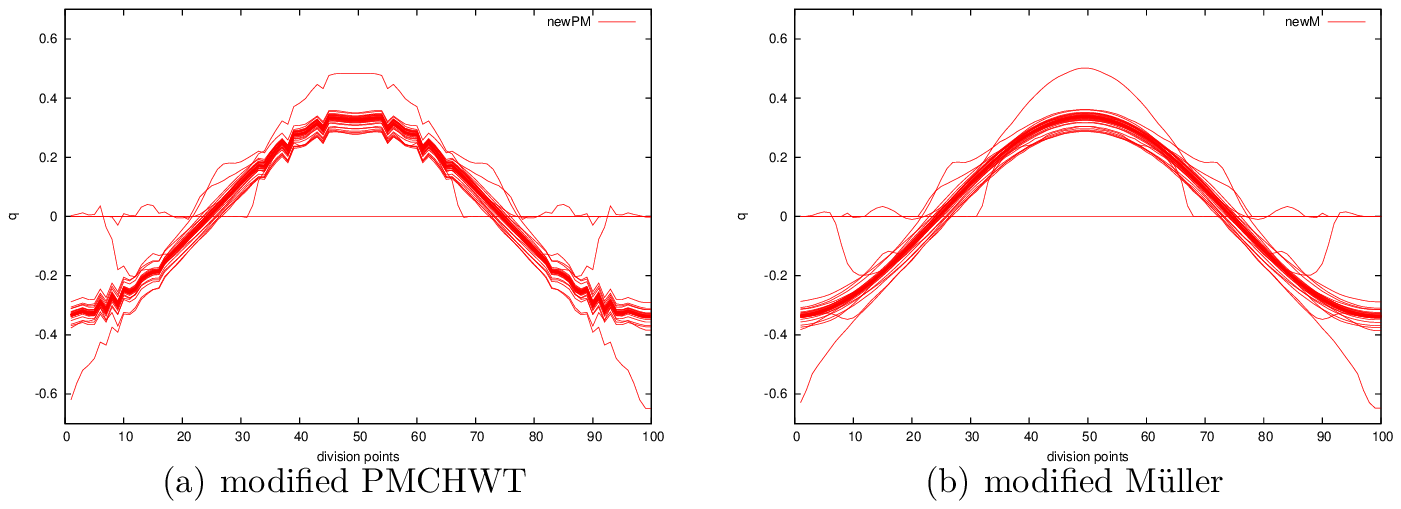}
    \end{center}
\caption{$q$ obtained with modified PMCHWT and M\"uller formulations}
  \label{lininc-transmission}
\end{figure}

To confirm this conclusion, we use the modified PMCHWT and M\"uller
equations to solve the same transmission problem as in Figs.\ref{old}
and \ref{new} after replacing the incident wave by the quasi-linear
one in \eqref{lininc}. As has been expected, the modified PMCHWT
result includes persistent noise, while the modified M\"uller result
is smooth as shown in Fig.\ref{lininc-transmission}. 
Further details of this subject will be presented elsewhere.

\section{Effects of space discretisation}
So far, we have neglected the effect of spatial discretisation in the
discussion of stability. This section discusses how the distribution
of the eigenvalues is influenced by the space discretisation.  We
restrict our attention to the circular scatterer case using the
Fourier series expansion in order to keep the discussion as analytical
as possible so that we can obtain insights.

\subsection{Formulation}
We consider a
circular boundary having a radius of $1$ with $N$ piecewise constant
arc elements whose endpoints (angles) are given as $\theta=\theta_1, \cdots, \theta_N$ ($\theta_{N+1}=\theta_{1}$). We consider an integral operator $K(\Omega)$ which maps a function $v$
on the boundary to another function on the boundary.
The function $v$ is then approximated with the piecewise constant basis function on
each element as follows:
\begin{align*}
  v(\theta)\approx \sum_{p=1}^N v_p N_p(\theta), \ N_p(\theta)=\left\{
    \begin{array}{cc}
      1 & \theta\in [\theta_p,\theta_{p+1}]\\
      0 & \mbox{otherwise}
    \end{array}
    \right.
\end{align*}
where $\theta$ is the  angular coordinate of the position on the boundary.  This is considered to be a reasonable approximation of the discretisation with straight line boundary
elements.
The basis function is now expanded into the Fourier
series given by:
\begin{align*}
  N_p(\theta)\approx\sum_{l=-M}^{M} V_l^p e^{il\psi}
\end{align*}
where  $M$ is the
truncation number of the infinite Fourier series and $V_l^p$ is the coefficient of
the Fourier series given as follows:
\begin{align*}
  V_l^p=\frac{1}{2\pi} \int_{0}^{2\pi} e^{-il\psi}N_p \,d\psi=\frac{1}{2\pi} \int_{\theta_p}^{\theta_{p+1}} e^{-il\psi} \,d\psi.
\end{align*}
Suppose that $K$ is an integral operator for layer potentials such as $\hat S$, $\hat M$, etc., in frequency domain. Then the value of $K v$
at a collocation
point $\Theta_n=(\theta_n+\theta_{n+1})/2$ is given as
 \begin{align}
   { K}v|_{\Theta_n,n=1,\cdots,N}=\sum_{p=1}^{N}\sum_{l=-M}^{M}
   e^{il\Theta_n} \{HJ\}_lV_l^p v_p.
 \end{align}
 where  $\{HJ\}_l$ stands for the product of Hankel and Bessel functions (or
 their derivatives) with
 appropriate coefficients. The $\{HJ\}$'s for the integral operators used in this paper are given as follows:
 \begin{align*}
  &  S^\nu\leftrightarrow H_n^{(1)}(k_\nu)J_n(k_\nu)\\
  &  \dot{S}^\nu\leftrightarrow -i\Omega H_n^{(1)}(k_\nu)J_n(k_\nu)\\
  &  D^{\nu-},\,D^{T\nu+}\leftrightarrow k_\nu H_n^{(1)'}(k_\nu)J_n(k_\nu)\\
  &  D^{\nu+},\,D^{T\nu-}\leftrightarrow k_\nu H_n^{(1)}(k_\nu)J'_n(k_\nu)\\
  &  M^\nu\leftrightarrow -k^2_\nu H_n^{(1)'}(k_\nu)J'_n(k_\nu)/i/\Omega
\end{align*}
The question of stability of the discretised integral operators in time domain
is now reduced to the nonlinear eigenvalue problem for $\Omega$
for the following matrix:
 \begin{align}
   {\mathsf U}\left(
   \begin{array}{ccc}
     \mathcal{D}_{-M}&&0\\
     &\ddots&\\
     0& &\mathcal{D}_M
   \end{array}
   \right){\mathsf V}\label{eq:op}
 \end{align}
 where $\mathsf V$ ($\mathsf U$) is a $(2M+1)\times N$ ($N \times (2M+1)$)
 matrix whose $(l,p)$ ($(n,l)$) components are given by 
 \begin{align*}
   {\mathsf V}_{lp}=V_l^p \quad ({\mathsf U}_{nl}=e^{il\Theta_n})
 \end{align*}
and
 \begin{align*}
\mathcal{D}_l= \sum_{m=-\infty}^{\infty} \{HJ\}_{l}\left(\Omega-\frac{2m\pi}{\Delta t}\right)\hat{\phi}_{\Delta t}\left(\Omega-\frac{2m\pi}{\Delta t}\right)
 \end{align*}
 The corresponding matrices for the boundary integral equations for transmission problems
 can be obtained similarly. We can now solve the non-linear eigenvalue problem
 given by
 \begin{align*}
   K(\Omega) v =0
 \end{align*}
 with the standard SSM.

 \subsection{Numerical experiments}

 We now show results of some numerical experiments. We 
 consider various integral operators on the unit circle, setting
 $\rho=1, \ s=1$ and $\Delta t=2\pi/100$, respectively.  The mesh on
 the boundary is uniform with $\theta_p=2\pi p/N, \ p=1,2,\cdots,N+1$.
 The number of boundary subdivision $N$ is set to be either $N=100$ or
 $N=200$.  Accordingly, the truncation number of the Fourier series
 $M$ is set to be $M=10N+N/2$.  Also, the infinite series
 $\sum_m \{HJ\}_m$ in \eqref{eq:op} and similar ones for transmission
 problems are truncated with 100 terms since results obtained with
 1000 terms were almost identical.  We note that it is not very easy
 to calculate Bessel and Hankel functions of high order with large
 arguments included in this calculation. This problem is handled with
 the help of Exflib, a well-known multiple-precision
 library\cite{exflib}.

\begin{figure}[htbp]
    \begin{center}
\includegraphics[width=\linewidth]{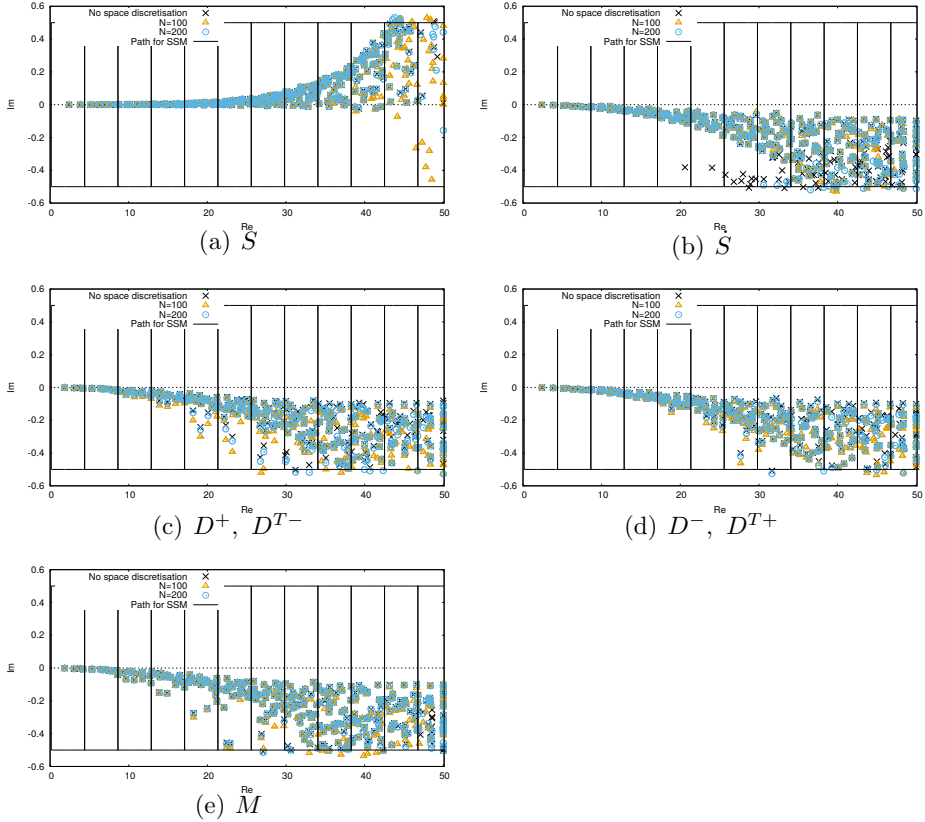}
      \caption[]{Eigenvalues of integral operators. $\times$: no
        space discretisation, $\triangle{}$: $N=100$, $\circ$: $N=200$}
    \label{fig:eig_op}
    \end{center}
  \end{figure}

 Figs.\,\ref{fig:eig_op} show the eigenvalue distributions of various
 integral operators considered in section \ref{stable_potentials}.
 The results with $N=100$ ($N=200$) are shown in triangular (circular)
 symbols and those without space discretisation (i.e., the
 characteristic roots given in previous sections, which we call ``no
 space discretisation'' in the rest of this paper) are given in cross
 symbols.  It is seen that the property of eigenvalue distributions
 does not change very much regardless of the space
 discretisations. Namely, stable potentials seem to remain stable for
 reasonable spatial divisions.  We also see that the $N=200$ results
 are closer to the no space discretisation results than those obtained
 with $N=100$.  These observations justify the use of the no space
 discretisation method in the discussion of the stability of the time
 domain BIEMs.  

\begin{figure}[htbp]
    \begin{center}
\includegraphics[width=\linewidth]{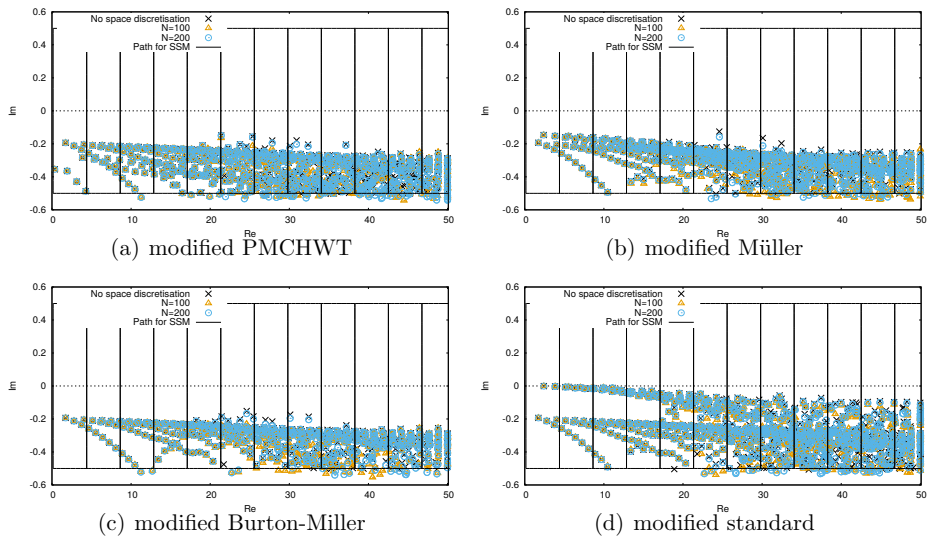}
      \caption[]{Eigenvalues of the modified integral equations for transmission
        problems. $\times$: no
        space discretisation, $\triangle{}$: $N=100$, $\circ$: $N=200$}
    \label{fig:eig_bie}
    \end{center}
  \end{figure}

  We next consider transmission problems.  Figs.\,\ref{fig:eig_bie}
  show the eigenvalues of the modified boundary integral equations for
  transmission problems, i.e., modified PMCHWT \eqref{PMnew}, M\"uller
  \eqref{Mnew}, BM \eqref{BMnew} and standard equations
  \eqref{snew}.  We set $\rho_1=1, \ \rho_2=0.37, \ s_1=1, \ s_2=0.2$,
  and $\Delta t=2\pi/100$, respectively as in \ref{ne_transmission}.
  Once again, the number of the spatial subdivision $N$ does not seem
  to change the distribution of the characteristic roots
  qualitatively.  We therefore conclude that the stability of these
  formulations can be inferred from the no space discretisation
  results.  Also, the finer the spatial discretisation the better
  approximation the no discretisation results become.  We thus expect
  that further spatial discretisation is not likely to affect the
  stability of the modified boundary integral equations.

\section{Non circular boundary}
\begin{figure}[htb]
    \begin{center}
\includegraphics[width=\linewidth]{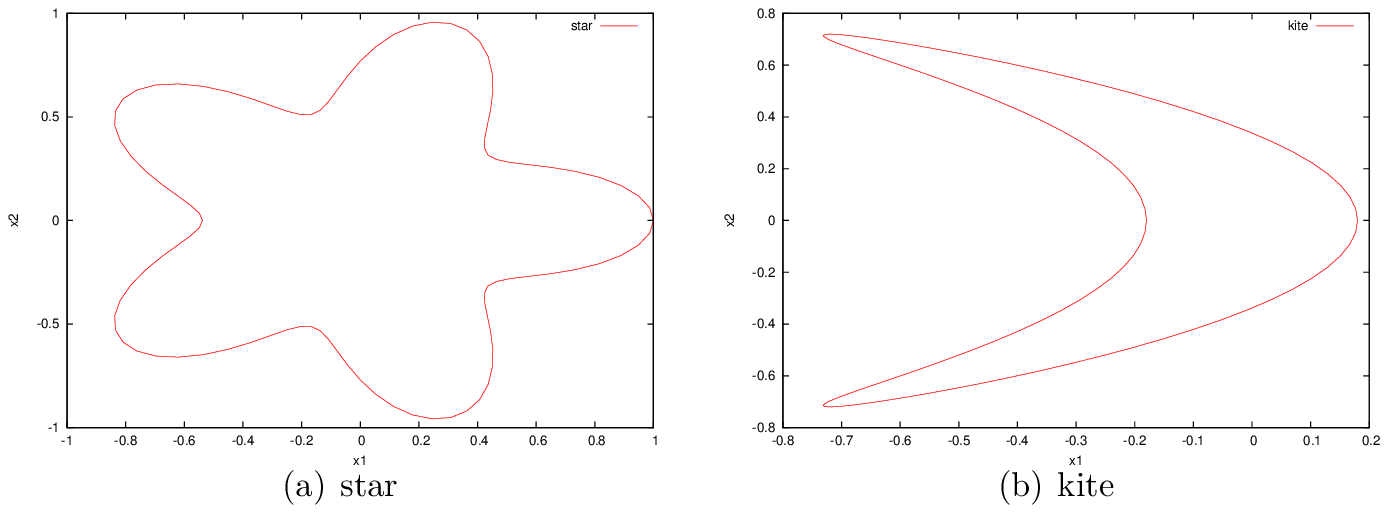}
      \caption[]{Domains}
    \label{fig:domains_sk}
    \end{center}
  \end{figure}

Finally, we test if the modified formulations remain stable for boundaries other than
circle. We consider transmission problems for a ``star'' (Fig.\ref{fig:domains_sk}(a)) given by
\begin{align*}
 (x_1,x_2)=((1+0.3\cos{5\theta})\cos{\theta}/1.3,(1+0.3\cos{5\theta})\sin{\theta}/1.3)
\end{align*}
and a ``kite''  (Fig.\ref{fig:domains_sk}(b)) given by
\begin{align*}
 (x_1,x_2)=(0.18(\cos{\theta}+2(\cos{2\theta}-1)),0.72\sin{\theta}) .
\end{align*}
The incident wave is the quadratic one in \eqref{quadinc} and
parameters such as material constants, number of boundary subdivision,
$\Delta t$ etc.\ are the same as those in the transmission problems
considered in the last section.  The boundary subdivision is uniform
with respect to $\theta$.

\begin{figure}[hbt]
    \begin{center}
\includegraphics[width=\linewidth]{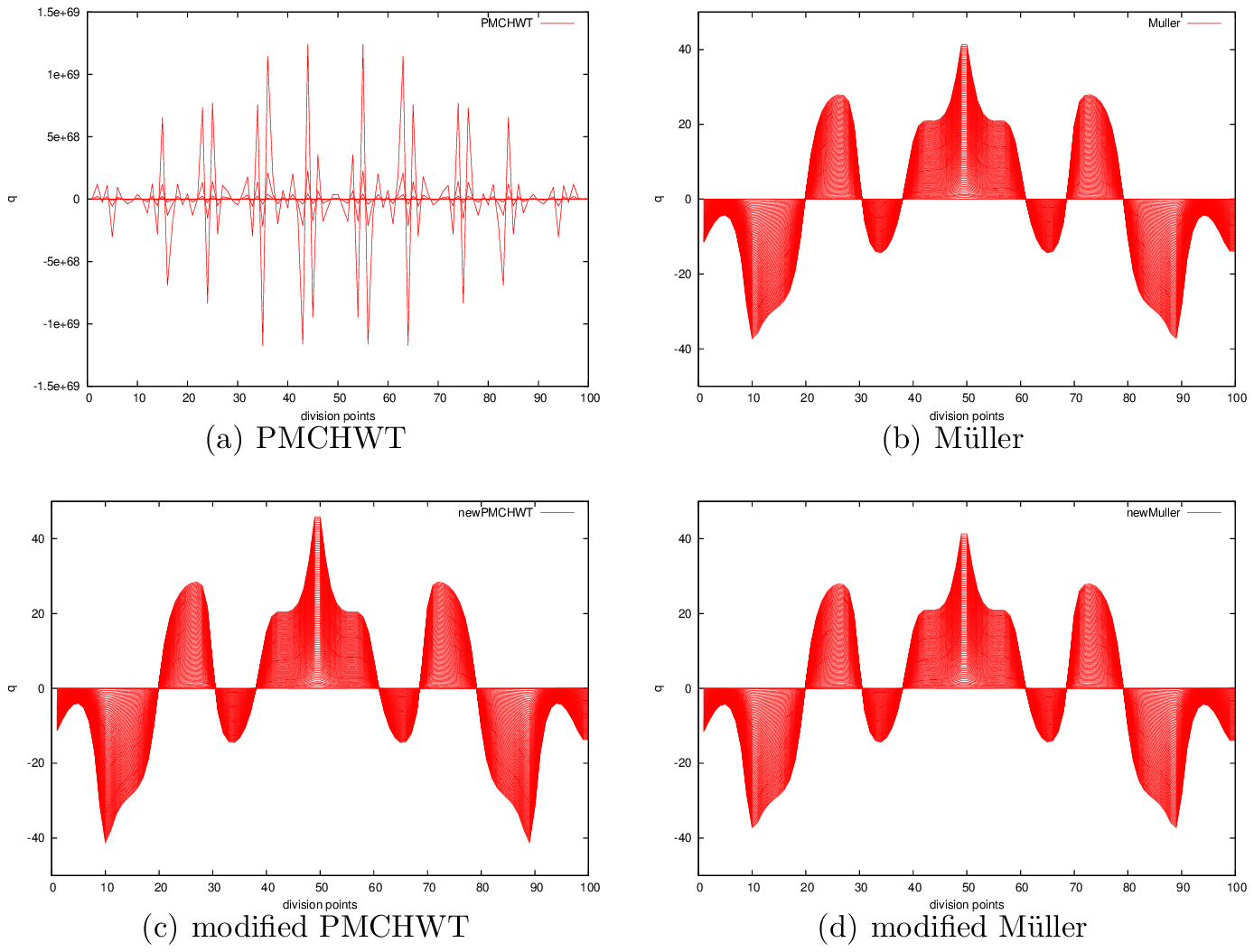}
      \caption[]{Transmission problems for ``star''}
    \label{fig:tr_st}
    \end{center}
  \end{figure}

\begin{figure}[htb]
    \begin{center}
\includegraphics[width=\linewidth]{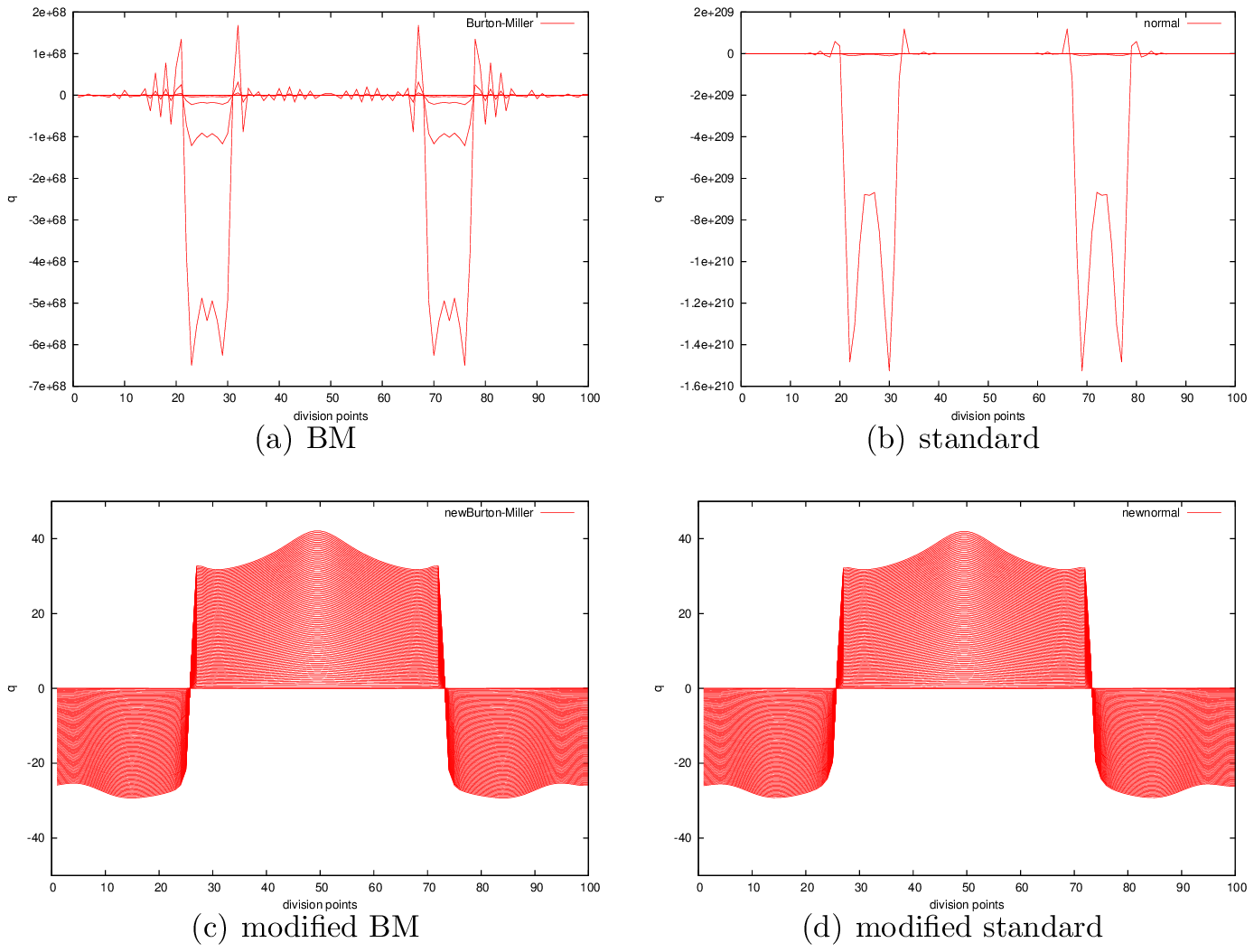}
      \caption[]{Transmission problems for ``kite''}
    \label{fig:tr_kt}
    \end{center}
  \end{figure}

  Fig.\ref{fig:tr_st} and Fig.\ref{fig:tr_kt} show the distribution of
  $q$ on the boundary obtained with various formulations.  Ordinary
  formulations except for M\"uller turn out to be unstable while
  modified formulations appear to be stable.  Also, the
  results obtained with modified formulations basically agree with
  each other except for details. These comments apply as well to other cases
which are not shown in the paper.

\section{Concluding remarks}

This paper revisited stability issues for BIEMs for the two
dimensional wave equation in time domain. We presented a stability
analysis based on integral equations in frequency domain and showed
its validity and usefulness in simple exterior or transmission
problems for circular domains.  The resulting non-linear eigenvalue
problems for the characteristic roots have been solved numerically
with SSM. We identified layer potentials which lead to stable integral
equations with linear time interpolation for our particular choices of
parameters. Combining these potentials, we could formulate stable
integral equations for transmission problems which include the
velocity and the normal flux of the solution on the boundary as
unknowns. Among integral equations considered the M\"uller formulation
was concluded to be a better choice in cases tested. All these
modified formulations were shown to remain stable in transmission
problems for star and kite shaped boundaries.

We remark that we have no intention to claim that the combination of
particular potentials always leads to stability or that the M\"uller
formulation is always the best choice in transmission problems. What
we showed in this paper is the fact that the proposed method of
stability analysis in frequency domain is useful in investigating the
stability and the accuracy of a time domain BIEM for the wave equation
in 2D, given particular integral equations and discretisation methods.
To be consistent with the purpose of this paper, we have restricted
our attention to simple problems where we can utilise analytical tools
as much as possible. Also, the numerical examples presented have been
limited to small number of cases. Two obvious next steps, therefore,
will be to consider more general boundaries and to carry out more
extensive numerical experiments. The former investigation will include
numerical treatment of \eqref{stab} in which one considers integral
equations having the function in \eqref{eq:ftg} as the kernel instead
of the fundamental solution in BIEs on general boundaries.  Another
interesting future direction is to test the modified formulations for
transmission problems in 3D. As a matter of fact, we have already
started investigations along this line. So far the same conclusions
for stability as have been presented here seem to hold in 3D as
well. Notice, however, that the use of the time domain stability
analysis based on \eqref{eq:dft} may be simpler than the frequency
domain approach in 3D (and, indeed, have already been utilised by many
authors including Walker et al.\cite{walker}, etc., as have been
mentioned) because of the finite ``tail'' of the fundamental
solution. We remark, however, that SSM will be useful in such
stability analysis in time domain as well because one may apply it
directly to a smaller eigenvalue problem in \eqref{algeh}. We are also
interested in the stability of interior problems in which true
eigenvalues may cause instability\cite{jang2}, thus requiring
different approaches than those utilised in this paper.  Finally, we
can mention investigations on the robustness of the algorithms as an
important future research subject.  As a matter of fact, we have
carefully tried to eliminate errors in the present investigation.  In
real world applications, however, one has to use numerical
integrations, truncated time steps, fast methods etc., which will
inevitably introduce errors.

\section*{Acknowledgement}
This work has been supported by JSPS KAKENHI Grant Number 18H03251.

\end{document}